# Karhunen-Loève expansions of mean-centered Wiener processes

**Paul Deheuvels[1]**

*L.S.T.A., Université Pierre et Marie Curie (Paris 6)*

**Abstract:** For $\gamma > -\frac{1}{2}$, we provide the Karhunen-Loève expansion of the weighted mean-centered Wiener process, defined by

$$W_\gamma(t) = \frac{1}{\sqrt{1+2\gamma}}\Big\{W\big(t^{1+2\gamma}\big) - \int_0^1 W\big(u^{1+2\gamma}\big)du\Big\},$$

for $t \in (0,1]$. We show that the orthogonal functions in these expansions have simple expressions in term of Bessel functions. Moreover, we obtain that the $L^2[0,1]$ norm of $W_\gamma$ is identical in distribution with the $L^2[0,1]$ norm of the weighted Brownian bridge $t^\gamma B(t)$.

## 1. Introduction and results

### 1.1. KL expansions of mean-centered Wiener processes

Let $\{W(t) : t \geq 0\}$ denote a standard Wiener process, and let

$$(1.1) \qquad \{B(t) : t \geq 0\} \stackrel{\text{law}}{=} \{W(t) - tW(1) : t \geq 0\},$$

denote a standard Brownian bridge, where "$\stackrel{\text{law}}{=}$" denotes equality in distribution. As a motivation to the present paper, we start by establishing, in Theorem 1.1 below, the Karhunen-Loève [KL] expansion of the *mean-centered Wiener process*

$$(1.2) \qquad W_0(t) := W(t) - \int_0^1 W(u)du, \quad t \in [0,1].$$

We recall the following basic facts about KL expansions (or representations). Let $d \geq 1$ be a positive integer. It is well-known (see, e.g., [1, 3, 9] and Ch.1 in [8]) that a centered Gaussian process $\{\zeta(\mathbf{t}) : \mathbf{t} \in [0,1]^d\}$, with covariance function $K_\zeta(\mathbf{s},\mathbf{t}) = E(\zeta(\mathbf{s})\zeta(\mathbf{t}))$ in $L^2\big([0,1]^d \times [0,1]^d\big)$, admits the (convergent in expected mean squares) KL expansion

$$(1.3) \qquad \zeta(t) \stackrel{\text{law}}{=} \sum_{k=1}^\infty \omega_k \sqrt{\lambda_k}\, e_k(\mathbf{t}),$$

where $\{\omega_k : k \geq 1\}$ is a sequence of independent and identically distributed [iid] normal $N(0,1)$ random variables, and the $\{e_k : k \geq 1\}$ form an orthonormal sequence

[1]L.S.T.A., Université Paris 6, 7 avenue du Château, F-92340 Bourg-la-Reine, France, e-mail: pd@ccr.jussieu.fr

in $L^2([0,1]^d)$, fulfilling

$$(1.4) \qquad \int_{[0,1]^d} e_k(\mathbf{t})e_\ell(\mathbf{t})d\mathbf{t} = \begin{cases} 1 & \text{if } k = \ell, \\ 0 & \text{if } k \neq \ell, \end{cases}$$

with $d\mathbf{t}$ denoting Lebesgue's measure. The $\lambda_1 \geq \lambda_2 \geq \cdots \geq 0$ in (1.3) are the eigenvalues of the Fredholm operator $h \in L^2([0,1]^d) \to T_\zeta h \in L^2([0,1]^d)$, defined via

$$(1.5) \qquad T_\zeta h(\mathbf{t}) = \int_{[0,1]^d} K_\zeta(\mathbf{s},\mathbf{t}) h(\mathbf{s}) d\mathbf{s} \quad \text{for} \quad t \in [0,1]^d.$$

We have, namely, $T_\zeta e_k = \lambda_k e_k$ for each $k \geq 1$. A natural consequence of the KL representation (1.3) is the distributional identity

$$(1.6) \qquad \int_{[0,1]^d} \zeta^2(\mathbf{t})d\mathbf{t} \stackrel{\text{law}}{=} \sum_{k=1}^{\infty} \lambda_k \omega_k^2.$$

Given these preliminaries, we may now state our first theorem.

**Theorem 1.1.** *The KL expansion of $\{W_0(t) : t \in [0,1]\}$ is given by*

$$(1.7) \qquad W_0(t) = W(t) - \int_0^1 W(u)du \stackrel{\text{law}}{=} \sum_{k=1}^{\infty} \omega_k \frac{\sqrt{2}\ \cos(k\pi t)}{k\pi} \quad \text{for} \quad t \in [0,1].$$

The proof of Theorem 1.1 is postponed until Section 2.1. In Remark 1.1 below, we will discuss some implications of this theorem in connection with well-known properties of the Brownian bridge, as defined in (1.1).

**Remark 1.1.** (1) We will provide, in the forthcoming Theorem 1.3 (in Section 1.4), a weighted version of the KL expansion of $W_0(\cdot)$, as given in (1.7). These KL expansions are new, up to our best knowledge.

(2) The well-known (see, e.g., [2] and [7]) KL expansion of the Brownian bridge (1.1) is very similar to (1.7), and given by

$$(1.8) \qquad B(t) = W(t) - tW(1) \stackrel{\text{law}}{=} \sum_{k=1}^{\infty} \omega_k \frac{\sqrt{2}\ \sin(k\pi t)}{k\pi} \quad \text{for} \quad 0 \leq t \leq 1.$$

As a direct consequence of (1.5)-(1.7)-(1.8), we obtain the distributional identities

$$(1.9) \qquad \int_0^1 \left\{W(t) - \int_0^1 W(u)du\right\}^2 dt \stackrel{\text{law}}{=} \int_0^1 B^2(t)dt \stackrel{\text{law}}{=} \sum_{k=1}^{\infty} \frac{\omega_k^2}{k^2\pi^2}.$$

The first identity in (1.9) is given, p.517 of [5], as a consequence of Fubini-Wiener arguments, in the spirit of the results of [6]. The second identity in (1.9) follows directly from (1.6), and is well-known (see, e.g., [4]).

The remainder of our paper is organized as follows. In Section 1.2 below, we mention, as a straightforward, but nevertheless useful, observation, that the knowledge of the distribution of the $L^2$ norm of a Gaussian process characterizes the eigenvalues of its KL expansion. In Section 1.3, we extend our study to Gaussian processes related to the Wiener sheet in $[0,1]^d$. The results of [5] will be instrumental in this case to establish a series of distributional identities of $L^2$ norms, between various Gaussian processes of interest. In Section 1.4, we provide the KL expansion of the *weighted mean-centered Wiener process* in dimension $d = 1$. The results of this section, are, as one could expect, closely related to [4], where related KL decompositions of weighted Wiener processes and Brownian bridges are established. In particular, in these results, we will make an instrumental use of Bessel functions. In the case of a general $d \geq 1$, we provide KL expansions for a general version of the mean-centered Wiener sheet. Finally, in Section 2, we will complete the proofs of the results given in Section 1.

### 1.2. The $L^2$ norm and KL eigenvalues

A question raised by (1.9) is as follows. Let $\zeta_1(\cdot)$ and $\zeta_2(\cdot)$ be two centered Gaussian processes on $[0,1]^d$, with covariance functions in $L^2\big([0,1]^d \times [0,1]^d\big)$, and KL expansions given by, for $\mathbf{t} \in [0,1]^d$,

$$(1.10) \qquad \zeta_j(\mathbf{t}) \stackrel{\text{law}}{=} \sum_{k=1}^{\infty} \omega_k \sqrt{\lambda_{k,j}}\, e_{k,j}(\mathbf{t}) \quad \text{for} \quad j = 1, 2.$$

Making use of the notation of [5], we write $\zeta_1 \stackrel{\text{Quad}}{=} \zeta_2$, when the $L^2([0,1]^d)$ norms of these two processes are identical in distribution. We may therefore write the equivalence

$$(1.11) \qquad \zeta_1 \stackrel{\text{Quad}}{=} \zeta_2 \quad \Leftrightarrow \quad \int_{[0,1]^d} \zeta_1^2(\mathbf{t}) d\mathbf{t} \stackrel{\text{law}}{=} \int_{[0,1]^d} \zeta_2^2(\mathbf{t}) d\mathbf{t}.$$

What can be said of the eigenvalue sequences $\{\lambda_{k,j} : k \geq 1\}$, $j = 1, 2$ when (1.11) is fulfilled? The answer to this question is given below.

**Theorem 1.2.** *The condition* $\zeta_1 \stackrel{\text{Quad}}{=} \zeta_2$ *is equivalent to the identity*

$$(1.12) \qquad \lambda_{k,1} = \lambda_{k,2} \quad \text{for all} \quad k \geq 1.$$

*Proof.* The fact that (1.12) implies $\zeta_1 \stackrel{\text{Quad}}{=} \zeta_2$ is trivial, via (1.6). A simple proof of the converse implication follows from the expression of the moment-generating functions [mgf] (see, e.g., pp. 60–61 in [4]), for $j = 1, 2$,

$$(1.13) \quad E\Big( \exp\Big\{ z \int_{[0,1]^d} \zeta_j^2(\mathbf{t}) d\mathbf{t} \Big\} \Big) = \prod_{k=1}^{\infty} \Big\{ \frac{1}{1 - 2z\lambda_{k,j}} \Big\}^{1/2} \quad \text{for} \quad \text{Re}(z) < \frac{1}{2\lambda_{1,j}}.$$

A variety of methods are available to show that the equality of the mgf's in (1.13), for $j = 1, 2$, implies (1.12). A simple argument, suggested by D. M. Mason (personal communication), uses the fact that the large deviations of the $L^2[0,1]$ norms of the $\zeta_j(\cdot)$ are governed by $\lambda_{1,j}$, $j = 1, 2$ (see, e.g., Lemma 1.1 in [4]). This allows to show that $\lambda_{1,1} = \lambda_{1,2}$. The proof is completed by a straightforward induction, after substracting from $\zeta_j(\cdot)$, $j = 1, 2$. the first components of their respective KL expansions. □

### 1.3. Multivariate Wiener sheets and Brownian Bridges

Inspired by Theorems 1.1-1.2 and Remark 1.1, we will devote the remainder of our paper to derive of a series of new KL expansions of Gaussian processes, in

*Karhunen-Loève Expansions* 65the spirit of (1.7). To give an additional motivation to the forthcoming results, we introduce the following notation and definitions. For any integer $d \geq 1$, the *d-variate Wiener process* (or *Brownian sheet*) is a centered Gaussian process $\mathbf{W}(\mathbf{t})$, defining a continuous function of $\mathbf{t} = (t_1, \ldots, t_d) \in \mathbb{R}_+^d$, with covariance function given by, for $\mathbf{s} = (s_1, \ldots, s_d) \in \mathbb{R}_+^d$ and $\mathbf{t} = (t_1, \ldots, t_d) \in \mathbb{R}_+^d$,

$$\text{(1.14)} \qquad E\big(\mathbf{W}(\mathbf{s})\mathbf{W}(\mathbf{t})\big) = \prod_{i=1}^{d} \big(s_i \wedge t_i\big).$$

For each $i = 1, \ldots, d$, denote, respectively, by $\Delta_i$, $\Sigma_i$ and $\Theta_i$, the operators which map a $L^1\big([0,1]^d\big)$ function $f(\mathbf{t})$ of $\mathbf{t} = (t_1, \ldots, t_d) \in [0,1]^d$, into, respectively,

$$\text{(1.15)} \qquad \Delta_i f(\mathbf{t}) = f(\mathbf{t}) - t_i f(t_1, \ldots, t_{i-1}, 1, t_{i+1}, \ldots, t_d),$$

$$\text{(1.16)} \qquad \Sigma_i f(\mathbf{t}) = f(\mathbf{t}) - \int_0^1 f(t_1, \ldots, t_i, \ldots, t_d) dt_i,$$

$$\text{(1.17)} \qquad \Theta_i f(\mathbf{t}) = f(\mathbf{t}) - f(t_1, \ldots, t_{i-1}, 1, t_{i+1}, \ldots, t_d).$$

The *d-variate tied-down Brownian bridge* is a centered Gaussian process $\mathbf{B}_*(\mathbf{t})$, which is a continuous function of $\mathbf{t} = (t_1, \ldots, t_d) \in [0,1]^d$, defined, in terms of the $d$-variate Wiener process $\mathbf{W}(\mathbf{t})$, via the distributional identity

$$\text{(1.18)} \qquad \mathbf{B}_*(\mathbf{t}) \stackrel{\text{law}}{=} \Delta_1 \circ \ldots \circ \Delta_d \mathbf{W}(\mathbf{t}),$$

where "$\circ$" denotes the composition of applications. We define likewise the *d-variate mean-centered Wiener sheet* by setting

$$\text{(1.19)} \qquad \mathbf{W}_M(\mathbf{t}) \stackrel{\text{law}}{=} \Sigma_1 \circ \ldots \circ \Sigma_d \mathbf{W}(\mathbf{t}).$$

In (1.18), we use the notation $\mathbf{B}_*$ to distinguish the tied-down Brownian bridge $\mathbf{B}_*$ from the usual *d-variate standard Brownian bridge* $\mathbf{B}$, the latter being classically defined via the distributional identity

$$\text{(1.20)} \qquad \mathbf{B}(\mathbf{t}) \stackrel{\text{law}}{=} \mathbf{W}(\mathbf{t}) - \Big\{\prod_{i=1}^{d} t_i\Big\} \mathbf{W}(1, \ldots, 1).$$

When $d = 1$ the processes $\mathbf{W}$ and $\mathbf{B}_*, \mathbf{B}$, reduce, respectively, to the standard Wiener process $W$, and to the standard Brownian bridge $B$, via the distributional identities

$$\text{(1.21)} \qquad \mathbf{W} \stackrel{\text{law}}{=} W \quad \text{and} \quad \mathbf{B}_* \stackrel{\text{law}}{=} \mathbf{B} \stackrel{\text{law}}{=} B.$$

On the other hand, for $d \geq 2$, the processes $\mathbf{B}_*$ and $\mathbf{B}$ have different distributions. In particular, their covariance functions are defined, respectively, for $\mathbf{s} = (s_1, \ldots, s_d) \in \mathbb{R}_+^d$ and $\mathbf{t} = (t_1, \ldots, t_d) \in \mathbb{R}_+^d$, by

$$\text{(1.22)} \qquad E\big(\mathbf{B}_*(\mathbf{s})\mathbf{B}_*(\mathbf{t})\big) = \prod_{i=1}^{d} \big\{s_i \wedge t_i - s_i t_i\big\},$$

$$\text{(1.23)} \qquad E\big(\mathbf{B}(\mathbf{s})\mathbf{B}(\mathbf{t})\big) = \prod_{i=1}^{d} \big\{s_i \wedge t_i\big\} - \prod_{i=1}^{d} \big\{s_i t_i\big\}.$$

In the sequel, we will use the boldface notation $\mathbf{W}$, $\mathbf{B}_*$ and $\mathbf{B}$ in the general $d$-variate framework when $d \geq 1$ is arbitrary, and limit the formerly used notation, $W$ and $B$, to the univariate case, when $d = 1$. Let now $\boldsymbol{\gamma} = (\gamma_1, \ldots, \gamma_d) \in \mathbb{R}^d \in (-1, \infty)$ be a vector of constants. We define, respectively, the *weighted Wiener sheet* $\mathbf{W}^{(\boldsymbol{\gamma})}$, and the *weighted tied-down Brownian bridge* $\mathbf{B}_*^{(\boldsymbol{\gamma})}$, by

$$(1.24) \quad \mathbf{W}^{(\boldsymbol{\gamma})}(\mathbf{t}) = t_1^{\gamma_1} \ldots t_d^{\gamma_d} \mathbf{W}(\mathbf{t}) \quad \text{and} \quad \mathbf{B}_*^{(\boldsymbol{\gamma})}(\mathbf{t}) = t_1^{\gamma_1} \ldots t_d^{\gamma_d} \mathbf{B}(\mathbf{t}),$$

when $\mathbf{t} = (t_1, \ldots, t_d) \in (0,1]^d$. For convenience, we set $\mathbf{W}^{(\boldsymbol{\gamma})}(\mathbf{t}) = \mathbf{B}_*^{(\boldsymbol{\gamma})}(\mathbf{t}) = 0$, when $t_i = 0$ for some $i = 1, \ldots, d$. Introduce now the *upper-tail Wiener sheet* defined, for $\mathbf{t} = (t_1, \ldots, t_d) \in [0,1]^d$, by

$$(1.25) \quad \widetilde{\mathbf{W}}^{(\boldsymbol{\gamma})}(\mathbf{t}) = \int_{[t_1,1] \times \cdots \times [t_d,1]} u_1^{\gamma_1} \ldots u_d^{\gamma_d} \mathbf{W}(du_1, \ldots, du_d).$$

It is readily checked that, whenever $\gamma_i > -\frac{1}{2}$ for $i = 1, \ldots, d$, we have, for all $\mathbf{s} = (s_1, \ldots, s_d) \in [0,1]^d$ and $\mathbf{t} = (t_1, \ldots, t_d) \in [0,1]^d$,

$$E\bigl(\widetilde{\mathbf{W}}^{(\boldsymbol{\gamma})}(\mathbf{s}) \widetilde{\mathbf{W}}^{(\boldsymbol{\gamma})}(\mathbf{t})\bigr) = \prod_{i=1}^{d} \int_{s_i \vee t_i}^{1} u_i^{2\gamma_i} du_i$$

$$= \prod_{i=1}^{d} \left\{ \frac{(1 - s_i^{2\gamma_i+1}) \wedge (1 - t_i^{2\gamma_i+1})}{1 + 2\gamma_i} \right\},$$

so that we have the distributional identity (see, e.g. (3.11), p.505 in [5])

$$(1.26) \quad \widetilde{\mathbf{W}}^{(\boldsymbol{\gamma})}(\mathbf{t}) \stackrel{\text{law}}{=} \left\{ \prod_{i=1}^{d} \frac{1}{1+2\gamma_i} \right\} \mathbf{W}\bigl(1 - t_1^{1+2\gamma_1}, \ldots, 1 - t_d^{1+2\gamma_d}\bigr).$$

The following additional distributional identities will be useful, in view of the definitions (1.16) and (1.17) of $\Sigma_i$ and $\Theta_i$, for $i = 1, \ldots, d$. We set, for convenience, $\mathbf{0} = (0, \ldots, 0)$.

We define the *mean-centered weighted upper-tail Wiener sheet* by setting, for $\mathbf{t} = (t_1, \ldots, t_d) \in [0,1]^d$,

$$(1.27) \quad \widetilde{\mathbf{W}}_M^{(\boldsymbol{\gamma})}(\mathbf{t}) = \Sigma_1 \circ \cdots \circ \Sigma_d \widetilde{\mathbf{W}}^{(\boldsymbol{\gamma})}(\mathbf{t})$$

**Lemma 1.1.** *For each $\mathbf{t} \in [0,1]^d$, we have*

$$(1.28) \quad \widetilde{\mathbf{W}}^{(\mathbf{0})}(\mathbf{t}) = \mathbf{W}(1 - t_1, \ldots, 1 - t_d) \stackrel{\text{law}}{=} \Theta_1 \circ \cdots \circ \Theta_d \mathbf{W}(\mathbf{t}),$$

*and*

$$(1.29) \quad \mathbf{W}_M \stackrel{\text{law}}{=} \Sigma_1 \circ \cdots \circ \Sigma_d \mathbf{W}(\mathbf{t}) \stackrel{\text{law}}{=} \Sigma_1 \circ \cdots \circ \Sigma_d \widetilde{\mathbf{W}}_M^{(\mathbf{0})} \stackrel{\text{law}}{=} \widetilde{\mathbf{W}}^{(\mathbf{0})}(\mathbf{t}).$$

*Proof.* The proofs of (1.28) and (1.29) are readily obtained by induction on $d$. For $d = 1$, we see that (1.28) is equivalent to the identity $W(1-t) \stackrel{\text{law}}{=} W(t) - W(1)$, which is obvious. Likewise, (1.29) reduces to the formula

$$W(t) - W(1) - \int_0^1 \{W(u) - W(1)\} du = W(t) - \int_0^1 W(u) du.$$

We now assume that (1.28) holds at rank $d-1$, so that

$$\mathbf{W}(1-t_1,\ldots,1-t_{d-1},t_d) \stackrel{\text{law}}{=} \Theta_1 \circ \ldots \circ \Theta_{d-1}\mathbf{W}(t_1,\ldots,t_d).$$

We now make use of the observation that

$$\mathbf{W}(1-t_1,\ldots,1-t_{d-1},1-t_d) \stackrel{\text{law}}{=} \Theta_d\mathbf{W}(1-t_1,\ldots,1-t_{d-1},t_d).$$

This, in combination with the easily verified fact that the operators $\Theta_1,\ldots,\Theta_d$ commute, readily establishes (1.28) at rank $d$. The completion of the proof of (1.29) is very similar and omitted. □

Under the above notation, Theorem 3.1 in [5] establishes the following fact.

**Fact 1.1.** *For $d=2$ and $\boldsymbol{\gamma} \in (-\frac{1}{2},\infty)^d$, we have*

$$(1.30) \qquad \mathbf{B}_*^{(\boldsymbol{\gamma})} \stackrel{\text{Quad}}{=} \widetilde{\mathbf{W}}_M^{(\boldsymbol{\gamma})}.$$

**Remark 1.2.** Theorem 3.1 in [5] establishes (1.30), and, likewise, Corollary 3.1 of [5], establishes the forthcoming (1.31), only for $d=2$. It is natural to extend the validity of (1.30) to the case of $d=1$. In the forthcoming Section 1.5, we will make use of KL expansions to establish a result of the kind.

**Remark 1.3.** For $d=2$, we may rewrite (1.30) into the distributional identity, for $\gamma > -\frac{1}{2}$ and $\delta > -\frac{1}{2}$,

$$(1.31) \quad \begin{aligned} &\int_0^1\!\!\int_0^1 s^{2\gamma}t^{2\delta}\Big\{\mathbf{W}(s,t) - s\mathbf{W}(1,t) - t\mathbf{W}(s,1) + st\mathbf{W}(1,1)\Big\}^2 \\ &\stackrel{\text{law}}{=} \int_0^1\!\!\int_0^1 \Big\{s^\gamma t^\delta \widetilde{\mathbf{W}}(s,t) - \int_0^1 u^\gamma t^\delta \widetilde{\mathbf{W}}(u,t)du \\ &\qquad - \int_0^1 s^\gamma v^\delta \widetilde{\mathbf{W}}(s,v)dv + \int_0^1\!\!\int_0^1 u^\gamma v^\delta \widetilde{\mathbf{W}}(u,v)dudv\Big\}^2 dsdt. \end{aligned}$$

By combining (1.29), with (1.31), taken with $\gamma = \delta = 0$, we obtain the distributional identity

$$(1.32) \quad \begin{aligned} &\int_0^1\!\!\int_0^1 \Big\{\mathbf{W}(s,t) - s\mathbf{W}(1,t) - t\mathbf{W}(s,1) + st\mathbf{W}(1,1)\Big\}^2 \\ &\stackrel{\text{law}}{=} \int_0^1\!\!\int_0^1 \Big\{\mathbf{W}(s,t) - \int_0^1 \mathbf{W}(u,t)du \\ &\qquad - \int_0^1 \mathbf{W}(s,v)dv + \int_0^1\!\!\int_0^1 \mathbf{W}(u,v)dudv\Big\}^2 dsdt. \end{aligned}$$

*1.4. Weighted mean-centered Wiener processes ($d=1$)*

In this section, we establish the KL expansion of the univariate *mean-centered weighted Wiener process*, defined, for $\gamma > -\frac{1}{2}$ and $t \in [0,1]$, by

$$(1.33) \qquad W_\gamma(t) = \frac{1}{\sqrt{1+2\gamma}}\Big\{W\big(t^{1+2\gamma}\big) - \int_0^1 W\big(u^{1+2\gamma}\big)du\Big\}.$$

In view of (1.24), when $\gamma = 0$, we have $W_\gamma \stackrel{\text{Quad}}{=} B^{(\gamma)}$, therefore, by Theorem 1.2, the eigenvalues of the KL expansions of $W_\gamma$ and $B^{(\gamma)}$ must coincide in this case. We will largely extend this result by giving, in Theorem 1.3 below, the complete KL expansion of $W_\gamma(\cdot)$ for an arbitrary $\gamma > -\frac{1}{2}$. We need first to recall a few basic facts about Bessel functions (we refer to [12], and to Section 2 in [4] for additional details and references).

For each $\nu \in \mathbb{R}$, the *Bessel function of the first kind*, of index $\nu$, is defined by

$$(1.34) \qquad J_\nu(x) = \left(\tfrac{1}{2}x\right)^\nu \sum_{k=0}^{\infty} \frac{\left(-\tfrac{1}{4}x^2\right)^k}{\Gamma(\nu+k+1)\Gamma(k+1)}.$$

To render this definition meaningful when $\nu \in \{-1, -2, \ldots\}$ is a negative integer, we use the convention that $a/\infty = 0$ for $a \neq 0$. Since, when $\nu = -n$ is a negative integer, $\Gamma(\nu+k+1) = \Gamma(n+k+1) = \infty$ for $k = 0, \ldots, n-1$, the corresponding terms in the series (1.34) vanish, allowing us to write

$$(1.35) \qquad J_{-n}(x) = (-1)^n J_n(x).$$

One of the most important properties of Bessel functions is related to the second order homogeneous differential equation

$$(1.36) \qquad x^2 y'' + xy' + (x^2 - \nu^2)y = 0.$$

When $\nu \in \mathbb{R}$ is noninteger, the Bessel functions $J_\nu$ and $J_{-\nu}$ provide a pair of linearly independent solutions of (1.36) on $(0, \infty)$. On the other hand, when $\nu = n$ is integer, $J_n(x)$ and $J_{-n}(x)$ are linearly dependent, via (1.35). To cover both of these cases, it is useful to introduce the *Bessel function of the second kind* of index $\nu$. Whenever $\nu$ is noninteger, this function is defined by

$$(1.37) \qquad Y_\nu(x) = \frac{J_\nu(x)\cos\nu\pi - J_{-\nu}(x)}{\sin\nu\pi},$$

and whenever $\nu = n$ is integer, we set

$$(1.38) \qquad Y_n(x) = \lim_{\nu \to n} \frac{J_\nu(x)\cos\nu\pi - J_{-\nu}(x)}{\sin\nu\pi}.$$

In view of the definitions (1.37)–(1.38), we see that, for an arbitrary $\nu \in \mathbb{R}$, $J_\nu$ and $Y_\nu$ provide a pair of linearly independent solutions of (1.36) on $(0, \infty)$. The behavior of the Bessel functions of first and second kind largely differ at 0. In particular, when $\nu > 0$, we have, as $x \downarrow 0$ (see, e.g., p.82 in [4]),

$$(1.39) \qquad J_\nu(x) = (1 + o(1))\frac{(\tfrac{1}{2})^\nu x^\nu}{\Gamma(\nu+1)},$$

and

$$(1.40) \qquad Y_\nu(x) = (1 + o(1))\frac{\Gamma(\nu)}{\pi}\left(\tfrac{1}{2}x\right)^{-\nu}.$$

When $\nu > -1$, the positive roots (or zeros) of $J_\nu$ are isolated (see, e.g., Fact 2.1 in [4]) and form an increasing sequence

$$(1.41) \qquad 0 < z_{\nu,1} < z_{\nu,2} < \cdots,$$

such that, as $k \to \infty$,

$$(1.42) \qquad z_{\nu,k} = \left\{k + \tfrac{1}{2}(\nu - \tfrac{1}{2})\right\}\pi + o(1).$$

The next fact is Theorem 1.4 of [4].

**Fact 1.2.** *For any $\gamma > -1$, or, equivalently, for each $\nu = 1/(2(1+\gamma)) > 0$, the KL expansion of $B^{(\gamma)}(t) = t^\gamma B(t)$ on $(0,1]$ is given by*

$$t^\gamma B(t) \stackrel{\text{law}}{=} \sum_{k=1}^\infty \omega_k \sqrt{\lambda_k}\, e_k(t), \tag{1.43}$$

*where the $\{\omega_k : k \geq 1\}$ are i.i.d. random variables, and, for $k \geq 1$ and $t \in (0,1]$,*

$$\lambda_k = \left\{\frac{2\nu}{z_{\nu,k}}\right\}^2 \quad \text{and} \quad e_k(t) = t^{\frac{1}{2\nu}-\frac{1}{2}}\left\{\frac{J_\nu\left(z_{\nu,k} t^{\frac{1}{2\nu}}\right)}{\sqrt{\nu}\, J_{\nu-1}(z_{\nu,k})}\right\}. \tag{1.44}$$

**Theorem 1.3.** *Let $\gamma > -\frac{1}{2}$, or, equivalently, let $0 < \nu = 1/(2(1+\gamma)) < 1$. Then, the KL expansion of $W_\gamma$ is given by*

$$W_\gamma(t) = \frac{1}{\sqrt{1+2\gamma}}\left\{W(t^{1+2\gamma}) - \int_0^1 W(u^{1+2\gamma})du\right\} \stackrel{\text{law}}{=} \sum_{k=1}^\infty \omega_k \sqrt{\lambda_k}\, e_k(t), \tag{1.45}$$

*where the $\{\omega_k : k \geq 1\}$ are i.i.d. random variables, and, for $k \geq 1$ and $t \in (0,1]$,*

$$\lambda_k = \left\{\frac{2\nu}{z_{\nu,k}}\right\}^2 \quad \text{and} \quad e_k(t) = at^{\frac{1}{2\nu}-\frac{1}{2}}\frac{J_{\nu-1}\left(z_{\nu,k}\, t^{\frac{1}{2\nu}}\right)}{\sqrt{\nu}J_{\nu-1}(z_{\nu,k})}. \tag{1.46}$$

The proof of Theorem 1.3 is given in Section 2.2. As an immediate consequence of this theorem, in combination with Fact 1.2, we obtain that the relation

$$W_\gamma \stackrel{\text{Quad}}{=} B^{(\gamma)}, \tag{1.47}$$

holds for each $\gamma > -\frac{1}{2}$. A simple change of variables transforms this formula into

$$\frac{1}{(1+2\gamma)^2}\int_0^1 t^{\frac{-2\gamma}{1+2\gamma}}\left\{W(t) - \frac{1}{2\gamma+1}\int_0^1 s^{\frac{-2\gamma}{1+2\gamma}}W(s)ds\right\}^2 dt \tag{1.48}$$

$$\stackrel{\text{law}}{=} \int_0^1 t^{2\gamma}\left\{W(t) - tW(1)\right\}^2 dt.$$

### 1.5. Mean-centered Wiener processes ($d \geq 1$)

In this section, we establish the KL expansion of the (unweighted) multivariate *mean-centered Wiener process*, defined, for $\mathbf{t} \in [0,1]^d$, by

$$\mathbf{W}_M(\mathbf{t}) \stackrel{\text{law}}{=} \Sigma_1 \circ \cdots \circ \Sigma_d \mathbf{W}(\mathbf{t}). \tag{1.49}$$

We obtain the following theorem.

**Theorem 1.4.** *The KL expansion of $\mathbf{W}_M$ is given by*

$$\mathbf{W}_M(\mathbf{t}) \stackrel{\text{law}}{=} \sum_{k_1=1}^\infty \cdots \sum_{k_d=1}^\infty \omega_{k_1,\ldots,k_d}\left\{\prod_{i=1}^d \left(\frac{\sqrt{2}\,\cos(k_i\pi t_i)}{k_i\pi}\right)\right\}, \tag{1.50}$$

*where $\{\omega_{k_1,\ldots,k_d} : k_1 \geq 0,\ldots,k_d \geq 0\}$ denotes an array of i.i.d. normal $N(0,1)$ random variables.*

The proof of Theorem 1.4 is postponed until Section 2.1.

## 2. Proofs

### 2.1. Proof of Theorems 1.1 and 1.4

In spite of the fact that Theorem 1.1 is a particular case of Theorem 1.3, it is useful to give details about its proof, to introduce the arguments which will be used later on, in the much more complex setup of weighted processes. We start with the following easy lemma. Below, we let $W_0(t)$ be as in (1.2).

**Lemma 2.1.** *The covariance function of $W_0(\cdot)$ is given, for $0 \le s, t \le 1$, by*

$$(2.1) \qquad K_{W_0}(s,t) = E(W_0(s)W_0(t)) = s \wedge t - s - t + \tfrac{1}{2}s^2 + \tfrac{1}{2}t^2 + \tfrac{1}{3}.$$

*Proof.* In view of (1.2), we have the chain of equalities, for $0 \le s, t \le 1$,

$$\begin{aligned}
E(W_0(s)W_0(t)) &= E\Big(W(s)W(t) - W(s)\int_0^1 W(u)du \\
&\quad - W(t)\int_0^1 W(u)du + \int_0^1\int_0^1 W(u)W(v)dudv\Big) \\
&= s \wedge t - \int_0^1 (s \wedge u)du - \int_0^1 (t \wedge u)du + \int_0^1\int_0^1 (u \wedge v)dudv,
\end{aligned}$$

from where (2.1) is obtained by elementary calculations. □

Recalling the definition (1.5) of $T_\zeta$, we set below $\zeta(\cdot) = W_0(\cdot)$.

**Lemma 2.2.** *Let $\{y(t) : 0 \le t \le 1\}$ denote an eigenfunction of the Fredholm transformation $T_{W_0}$, pertaining to the eigenvalue $\lambda > 0$. Then, $y(\cdot)$ is infinitely differentiable on $[0,1]$, and a solution of the differential equation*

$$(2.2) \qquad \lambda y''(t) + y(t) = \int_0^1 y(u)du,$$

*subject to the boundary conditions*

$$(2.3) \qquad y'(0) = y'(1) = 0.$$

*Proof.* By (2.1), we have, for each $t \in [0,1]$,

$$(2.4) \qquad \begin{aligned}
\lambda y(t) &= \int_0^t sy(s)ds + t\int_t^1 y(s)ds + (\tfrac{1}{2}t^2 - t)\int_0^1 y(s)ds \\
&\quad + \int_0^1 (\tfrac{1}{2}s^2 - s + \tfrac{1}{3})y(s)ds.
\end{aligned}$$

It is readily checked that the RHS of (2.4) is a continuous function of $t \in [0,1]$. This, together with the condition that $\lambda > 0$ entails, in turn, that $y(\cdot)$ is continuous on $[0,1]$. By repeating this argument, a straightforward induction implies that $y(\cdot)$ is infinitely differentiable on $[0,1]$. This allows us to derivate both sides of (2.4) with respect to $t$. We so obtain that, for $t \in [0,1]$,

$$(2.5) \qquad \lambda y'(t) = (t-1)\int_0^1 y(s)ds + \int_t^1 y(s)ds.$$

By setting, successively, $t = 0$ and $t = 1$ in (2.5), we get (2.3). Finally, (2.2) is obtained by derivating both sides of (2.5) with respect to $t$. □

*Proof of Theorem 1.1.* It is readily checked that the general solution of the differential equation (2.2) is of the form

$$y(t) = a \cos\left(\frac{t}{\sqrt{\lambda}} + b\right) + c, \tag{2.6}$$

for arbitrary choices of the constants $a, b, c \in \mathbb{R}$. The limit conditions (2.3) imply that, in (2.6), we must restrict $b$ and $\lambda$ to fulfill $b = 0$ and $\lambda \in \{1/(k^2\pi^2) : k \geq 1\}$. We now check that the function $y(t) = a\cos(k\pi t)$ is a solution of the equation (2.4), taken with $\lambda = 1/(k^2\pi^2)$. Towards this aim, we first notice that $y(t) = a\cos(k\pi t)$ fulfills

$$\int_0^1 y(s)ds = \frac{a}{k\pi}\Big[\sin(k\pi t)\Big]_{t=0}^{t=1} = 0. \tag{2.7}$$

This, in turn, readily entails that $y(t) = a\cos(k\pi t)$ satisfies the relation (2.5). By integrating both sides of this relation with respect to $t$, we obtain, in turn, that $y(t) = a\cos(k\pi t)$ fulfills

$$\lambda y(t) = \int_0^t sy(s)ds + t\int_t^1 y(s)ds + (\tfrac{1}{2}t^2 - t)\int_0^1 y(s)ds$$
$$+ \int_0^1 (\tfrac{1}{2}s^2 - s + \tfrac{1}{3})y(s)ds + C,$$

for some constant $C \in \mathbb{R}$. All we need is to check that $C = 0$ for some particular value of $t$. If we set $t = 1$, and make use of (2.7), this reduces to show, by integration by parts, that

$$\frac{\cos(k\pi)}{k^2\pi^2} = \frac{1}{2}\int_0^1 s^2 y(s)ds = \frac{1}{2k\pi}\Big[t^2\sin(k\pi t)\Big]_{t=0}^{t=1} - \frac{1}{k\pi}\int_0^1 s\sin(k\pi s)ds$$
$$= \frac{1}{(k\pi)^2}\Big[t\cos(k\pi t)\Big]_{t=0}^{t=1} - \frac{1}{(k\pi)^2}\int_0^1 \cos(k\pi s)ds = \frac{\cos(k\pi)}{k^2\pi^2}.$$

The just-proved fact that $y(t) = a\cos(k\pi t)$ satisfies the relation (2.5) implies, in turn, that this same equation, taken with $y(t) = a\cos(k\pi t) + c$, reduces to

$$\frac{c}{k^2\pi^2} = c\int_0^1 (\tfrac{1}{2}s^2 - s + \tfrac{1}{3})ds = 0.$$

Since this relation is only possible for $c = 0$, we conclude that the eigenfunctions of $T_{W_0}$ are of the form $y(t) = a\cos(k\pi t)$, for an arbitrary $a \in \mathbb{R}$. Given this fact, the remainder of the proof of Theorem 1.1 is straightforward. □

To establish Theorem 1.4, we will make use of the following lemma. We let $\mathbf{W}_M$ be defined as in (1.19), and recall (2.1).

**Lemma 2.3.** *We have, for* $\mathbf{s} = (s_1, \ldots, s_d) \in [0,1]^d$ *and* $\mathbf{t} = (t_1, \ldots, t_d) \in [0,1]^d$,

$$E\big(\mathbf{W}_M(\mathbf{s})\mathbf{W}_M(\mathbf{s})\big) = \prod_{i=1}^d K_{W_0}(s_i, t_i). \tag{2.8}$$

*Proof.* It follows from the following simple argument. We will show that, for an arbitrary $1 \leq j \leq d$,

$$
\begin{aligned}
&E\big(\{\Sigma_1 \circ \ldots \circ \Sigma_j \mathbf{W}(\mathbf{s})\}\{\Sigma_1 \circ \ldots \circ \Sigma_j \mathbf{W}(\mathbf{t})\}\big) \\
&= \Big\{\prod_{i=1}^{j} K_{W_0}(s_i, t_i)\Big\} \prod_{i=j+1}^{d} \{s_i \wedge t_i\}.
\end{aligned}
\tag{2.9}
$$

Since $\Sigma_1, \ldots, \Sigma_d$ are linear mappings, the proof of (2.9) is readily obtained by induction on $j = 1, \ldots, d$. This, in turn, implies (2.8) for $j = d$. □

*Proof of Theorem 1.4.* In view of (2.8), the proof of the theorm follows readily from a repeated use of Lemma 4.1 in [5], which we state below for convenience. Let $\zeta_1(\mathbf{s})$, $\zeta_2(\mathbf{t})$, and $\zeta_3(\mathbf{s}, \mathbf{t})$ be three centered Gaussian processes, functions of $\mathbf{s} \in [0,1]^p$ and $\mathbf{t} \in [0,1]^q$. We assume that

$$E\big(\zeta_3(\mathbf{s}', \mathbf{t}'))\zeta_3(\mathbf{s}'', \mathbf{t}'')\big) = E\big(\zeta_1(\mathbf{s}'))\zeta_1(\mathbf{s}'')\big) E\big(\zeta_2(\mathbf{t}'))\zeta_2(\mathbf{t}'')\big).$$

Then, if the KL expansions of $\zeta_1$ and $\zeta_2$ are given by

$$\zeta_1(\mathbf{s}) \stackrel{\text{law}}{=} \sum_{k=1}^{\infty} \omega_k \sqrt{\lambda_{k,1}}\, e_{k,1}(\mathbf{s}) \quad \text{and} \quad \zeta_2(\mathbf{t}) \stackrel{\text{law}}{=} \sum_{k=1}^{\infty} \omega_k \sqrt{\lambda_{k,2}}\, e_{k,2}(\mathbf{t}),$$

the KL expansion of $\zeta_3$ is given by

$$\zeta_3(\mathbf{s}, \mathbf{t}) \stackrel{\text{law}}{=} \sum_{k=1}^{\infty} \sum_{\ell=1}^{\infty} \omega_{k,\ell} \sqrt{\lambda_{k,1}\lambda_{k,2}}\, e_{k,1}(\mathbf{s})e_{k,2}(\mathbf{t}),$$

where $\{\omega_{k,\ell} : k \geq 1, \ell \geq 1\}$ denotes an array of iid normal $N(0,1)$ r.v.'s. □

### 2.2. Proof of Theorem 1.3

Recall that

$$W_\gamma(t) = \frac{1}{\sqrt{1+2\gamma}}\Big\{ W(t^{1+2\gamma}) - \int_0^1 W(u^{1+2\gamma})du \Big\}. \tag{2.10}$$

**Lemma 2.4.** *The covariance function of $W_\gamma(\cdot)$ is given, for $0 \leq s, t \leq 1$, by*

$$
\begin{aligned}
K_{W_\gamma}(s,t) &= E(W_\gamma(s)W_\gamma(t)) \\
&= \frac{1}{2\gamma+1}\Big\{(s \wedge t)^{1+2\gamma} - s^{1+2\gamma} - t^{1+2\gamma} \\
&\quad + \frac{1+2\gamma}{2+2\gamma} s^{2+2\gamma} + \frac{1+2\gamma}{2+2\gamma} t^{2+2\gamma} + \frac{2}{(2+2\gamma)(3+2\gamma)}\Big\}.
\end{aligned}
\tag{2.11}
$$

*Proof.* The proof of (2.11) being very similar to the above given proof of (2.1), we omit details. □

**Lemma 2.5.** *Fix any $\gamma > -\frac{1}{2}$. Let $e(t)$ be an eigenfunction of $T_{W_\gamma}$ pertaining to a positive eigenvalue, $\lambda > 0$, of this operator. Then, the function $y(x) = e(x^{1/(1+2\gamma)})$ is solution on $[0,1]$ of the differential equation*

$$(1+2\gamma)^2 \lambda x^{\frac{2\gamma}{1+2\gamma}}\, y''(x) + y(x) = \frac{1}{1+2\gamma}\int_0^1 u^{-\frac{2\gamma}{1+2\gamma}}\, y(u)du, \tag{2.12}$$

*with boundary conditions*

(2.13) $$y'(0) = y'(1) = 0.$$

*Proof.* Let $e(t)$ fulfill, for some $\lambda > 0$,

$$\lambda e(t) = \int_0^1 K_{W_\gamma}(s,t)e(s)ds.$$

By (2.11), we may rewrite this identity into

(2.14)
$$\begin{aligned}(2\gamma+1)\lambda e(t) &= \int_0^t s^{1+2\gamma}e(s)ds + t^{1+2\gamma}\int_t^1 e(s)ds \\ &- \int_0^1 s^{1+2\gamma}e(s)ds - t^{1+2\gamma}\int_0^1 e(s)ds + \left\{\frac{1+2\gamma}{2+2\gamma}\right\}\int_0^1 s^{2+2\gamma}e(s)ds \\ &+ \left\{\frac{1+2\gamma}{2+2\gamma}\right\}t^{2+2\gamma}\int_0^1 e(s)ds + \frac{1}{(1+\gamma)(3+2\gamma)}\int_0^1 e(s)ds.\end{aligned}$$

A straightforward induction shows readily that any function $e(\cdot)$ fulfilling (2.14) is infinitely differentiable on $[0,1]$. This, in turn, allows us to derivate both sides of (2.14), as to obtain the equation

(2.15) $$\lambda e'(t) = -t^{2\gamma}\int_0^t e(s)ds + t^{1+2\gamma}\int_0^1 e(s)ds.$$

Set now $e(t) = y(t^{1+2\gamma})$ in (2.15). By changing variables in the LHS of (2.15), we obtain the equation

(2.16) $$(1+2\gamma)\lambda y'(t^{1+2\gamma}) = -\int_0^t e(s)ds + t\int_0^1 e(s)ds.$$

The relation (2.13) follows readily from (2.16), taken, successively, for $t=0$ and $t=1$. By derivating both sides of (2.16), and after setting $e(s) = y(s^{1+2\gamma}))$, we get

$$(1+2\gamma)^2\lambda y''(t^{1+2\gamma})t^{2\gamma} + y(t^{1+2\gamma}) = \int_0^1 y(s^{1+2\gamma})ds.$$

After making the change of variables $t = x^{1/(1+2\gamma)}$ and $s = u^{1/(1+2\gamma)}$, we may rewrite this last equation into

$$(1+2\gamma)^2\lambda x^{\frac{2\gamma}{1+2\gamma}}\,y''(x) + y(x) = \frac{1}{1+2\gamma}\int_0^1 u^{-\frac{2\gamma}{1+2\gamma}}\,y(u)du,$$

which is (2.16). The proof of Lemma 2.5 is now completed. □

Recall the definitions (1.34) and (1.37)–(1.38). In view of (2.12) the following fact will be instrumental for our needs (refer to Fact 2.3 in [4], and p.666 in [10]).

**Fact 2.1.** *Let $\lambda > 0$ and $\beta > -1$ be real constants. Then, the differential equation*

(2.17) $$\lambda y''(x) + x^{2\beta}y(x) = 0,$$

*has fundamental solutions on $(0,\infty)$ given by*

(2.18) $$x^{1/2}J_{\frac{1}{2(\beta+1)}}\left(\frac{x^{\beta+1}}{(\beta+1)\sqrt{\lambda}}\right) \quad \text{and} \quad x^{1/2}Y_{\frac{1}{2(\beta+1)}}\left(\frac{x^{\beta+1}}{(\beta+1)\sqrt{\lambda}}\right),$$

when $1/(2(\beta+1)) \in \mathbb{R}$, and

$$(2.19) \quad x^{1/2} J_{\frac{1}{2(\beta+1)}}\left(\frac{x^{\beta+1}}{(\beta+1)\sqrt{\lambda}}\right) \quad \text{and} \quad x^{1/2} J_{-\frac{1}{2(\beta+1)}}\left(\frac{x^{\beta+1}}{(\beta+1)\sqrt{\lambda}}\right),$$

when $1/(2(\beta+1))$ is noninteger.

**Lemma 2.6.** *Assume that $\gamma > -\frac{1}{2}$. Then, the solutions on $[0,1]$ of the differential equation*

$$(2.20) \quad (1+2\gamma)^2 \lambda x^{\frac{2\gamma}{1+2\gamma}} y''(x) + y(x) = 0,$$

*with boundary conditions*

$$(2.21) \quad y'(0) = y'(1) = 0,$$

*are of the form*

$$(2.22) \quad y(x) = a x^{1/2} J_{-\frac{1+2\gamma}{2(1+\gamma)}}\left(\frac{x^{\frac{1+\gamma}{1+2\gamma}}}{(1+\gamma)\sqrt{\lambda}}\right),$$

*for some arbitrary constant $a \in \mathbb{R}$.*

*Proof.* Set $\beta = -\gamma/(1+2\gamma)$ in Fact 2.1, and observe that the assumption that $\gamma > -\frac{1}{2}$ implies that $\beta > -1$, and that $1/(2(\beta+1))$ is noninteger. By Fact 2.1, the general solution on $(0,1]$ of the homogeneous differential equation (2.20) is of the form

$$(2.23) \quad y(x) = b x^{1/2} J_{\frac{1+2\gamma}{2(1+\gamma)}}\left(\frac{x^{\frac{1+\gamma}{1+2\gamma}}}{(1+\gamma)\sqrt{\lambda}}\right) + a x^{1/2} J_{-\frac{1+2\gamma}{2(1+\gamma)}}\left(\frac{x^{\frac{1+\gamma}{1+2\gamma}}}{(1+\gamma)\sqrt{\lambda}}\right),$$

where $a$ and $b$ are arbitrary constants. It is straightforward, given (1.39) and (1.40), that, for some constant $\rho_1 > 0$, as $x \downarrow 0$,

$$(2.24) \quad x^{1/2} J_{\frac{1+2\gamma}{2(1+\gamma)}}\left(\frac{x^{\frac{1+\gamma}{1+2\gamma}}}{(1+\gamma)\sqrt{\lambda}}\right) = (1+o(1))\rho_1 x^{\frac{1}{2}+\left\{\frac{1+2\gamma}{2(1+\gamma)}\right\}\frac{1+\gamma}{1+2\gamma}} \sim \rho_1 x,$$

whereas, for some constant $\rho_2 > 0$,

$$(2.25) \quad x^{1/2} Y_{-\frac{1+2\gamma}{2(1+\gamma)}}\left(\frac{x^{\frac{1+\gamma}{1+2\gamma}}}{(1+\gamma)\sqrt{\lambda}}\right) = (1+o(1))\rho_2 x^{\frac{1}{2}-\left\{\frac{1+2\gamma}{2(1+\gamma)}\right\}\frac{1+\gamma}{1+2\gamma}} \to \rho_2.$$

In view of (2.21), these relations imply that $b = 0$, so that (2.22) now follows from (2.23). $\square$

*Proof of Theorem 1.3.* Set $\nu = 1/(2(1+\gamma))$. Obviously, the condition that $\gamma > -\frac{1}{2}$ is equivalent to $0 < \nu < 1$. We will show that the eigenvalues of $T_{W_\gamma}$ are given by

$$(2.26) \quad \lambda_k = \left\{\frac{2\nu}{z_{\nu,k}}\right\} \quad \text{for} \quad k = 1, 2 \ldots .$$

Given this notation, it is readily checked that

$$\frac{1+2\gamma}{2(1+\gamma)} = 1 - \nu \quad \text{and} \quad 1 + 2\gamma = \frac{1}{\nu} - 1.$$

Therefore, letting $y(x)$ be as in (2.22), and setting $x = t^{1+2\gamma}$, we see that the eigenfunction $e_k(t)$ of $T_{W_\gamma}$ pertaining to $\lambda_k$ is of the form

$$(2.27) \qquad e_k(t) = e(t) = y(t^{1+2\gamma}) = at^{\frac{1}{2\nu}-\frac{1}{2}} J_{\nu-1}\left(z_{\nu,k}\, t^{\frac{1}{2\nu}}\right),$$

for some $a \in \mathbb{R}$. Here, we may use Formula 50:10:2, p.529 in [11], namely, the relation

$$\frac{d}{dx}\left\{x^\rho J_\rho(x)\right\} = x^\rho J_{\rho+1}(x),$$

to show that, for some appropriate constant $C$,

$$(2.28) \qquad \frac{d}{dx}\left\{x^{1-\nu} J_{\nu-1}(z_{\nu,k}\, x)\right\} = Cx^{1-\nu} J_\nu(z_{\nu,k}\, x) = 0,$$

for $x = 1$ and $x = 0$. This, in turn, allows to check that the function $y(x)$ in (2.27) fulfills (2.21). Now, we infer from (2.27), in combination with Fact 2.2, p.84 in [4], that

$$1 = \int_0^1 e_k^2(t)dt = 2\nu a^2 \int_0^1 t^{\frac{1}{2\nu}} J_{\nu-1}^2\left(z_{\nu,k}\, t^{\frac{1}{2\nu}}\right) dt^{\frac{1}{2\nu}}$$
$$= 2\nu a^2 \int_0^1 u J_{\nu-1}^2(z_{\nu,k}\, u)du = \nu a^2 J_{\nu-1}^2(z_{\nu,k}).$$

Given this last result, the completion of the proof of Theorem 1.3 is straightforward. □